\documentclass[a4paper,12pt]{amsart}
\newtheorem{thm}{Theorem}

\newtheorem{prop}[thm]{Proposition}
\newtheorem{lemme}[thm]{Lemma}

\usepackage{amssymb}
\usepackage{amsmath}
\usepackage[all]{xy}
\usepackage{a4wide}
\usepackage[dvips]{graphicx}
\newcommand{\Cb}{\mathbb{C}}
\newcommand{\Pb}{\mathbb{P}}

\newcommand{\Ar}{\mathcal{A}}
\newcommand{\Br}{\mathcal{B}}

\newcommand{\Zb}{\mathbb{Z}}

\newcommand{\Tb}{\mathbb{T}}

\newcommand{\Lr}{\mathcal{L}}

\title{Non-formality of Milnor fibers of line arrangements}
\author{Hugues Zuber}
\date{}
\address{Laboratoire J.A. Dieudonn\'e, UMR du CNRS 6621, Universit\'e de Nice - Sophia Antipolis, Parc Valrose, 06108 Nice Cedex 02, France}
\email{zuber@unice.fr}
\begin{document}
\begin{abstract}
The complement $M$ of a complex line arrangement $\Ar$ in $\Pb^2(\Cb)$ is formal, and admits a finite cover $F \rightarrow M$, where $F$ is the associated global Milnor fiber. An example where $F$ is not formal (not even $1$-formal), is given in this article.
\end{abstract}
\maketitle

\section{Introduction}

Let $f \in \Cb[x_0,\dots, x_n]$ be a homogeneous polynomial. The natural projection $\pi : \Cb^{n+1}\setminus \{0\} \rightarrow \Pb^n(\Cb)$ restricts to a map $\pi : F \longrightarrow M$, from the global Milnor fiber $F = f^{-1}(1)$ to the complement $M = \Pb^n(\Cb) \setminus f^{-1}(0)$. Consider the particular case when $f$ splits into a product of linear forms, so that $f^{-1}(0) \subset \Pb^n$ is a hyperplane arrangement. Then it is known that $M$ is formal \cite{Br73}, and one may wonder if this property passes to the Milnor fiber $F$ (see Question 2.10 in \cite{DP09}). The goal of this paper is to give an example where $F$ is not formal (and actually, not even $1$-formal).\\

The characteristic varieties (of rank one local systems) of a complex variety $X$ are defined as the jumping loci for the cohomology of complex rank one local systems on $X$:
$$ V_1^d(X)) = \{ \Lr \in H^1(X,\Cb^{\ast}),\, dim_{\Cb} H^1(X,\Lr) \ge d \}. $$
These are algebraic subvarieties of the algebraic group
$$ \Tb(X) = H^1(X,\Cb^{\ast}) \simeq Hom(H_1(X),\Cb^{\ast}) $$
whose connected component of identity is $\Tb(X)_1 \simeq (\Cb^{\ast})^{b_1(X)}$, where $b_1(X)$ is the first Betti number of $X$.\\
Another class of invariants, the resonance varieties, are defined at the level of $H^1(X,\Cb)$ by
$$ R_1^d(X) = \{ \omega \in H^1(X,\Cb),\, dim_{\Cb} H^1(H^{\bullet}(X,\Cb),\omega\wedge) \ge d \}. $$
Actually, in many cases, those two kinds of varieties are related via the exponential map
$$ exp : H^1(M,\Cb) \longrightarrow H^1(M,\Cb^{\ast}).$$
Precisely, if the complex variety $X$ is $1$-formal, then $exp(R_1^d(X))$ is the reunion of the components of $V_1^d(X)$ that contain the origin $1 \in H^1(X,\Cb^{\ast})$. So, in this case, $R_1^d(X)$ is the tangent cone at the origin of $V_1^d(X)$, $TC_1V_1^d(X)$. This theorem was proved by A. Dimca, S. Papadima, and A. Suciu in \cite{DPS08}.\\
The example given in this paper is a line arrangement in $\Pb^2(\Cb)$ such that $R_1^1(F) \neq TC_1V_1^1(F)$. It is defined by $(f=0)$, with
$$ f = (x^3-y^3)(x^3-z^3)(y^3-z^3). $$
It is described in \cite{BDS09}, 3.2 (i), and also in \cite{D92}, Theorem 6.4.15. This arrangement has nine lines, no ordinary double points, and twelve triple points. The geometric monodromy on the Milnor fiber $h: F \rightarrow F$ acts on $H^1(F,\Cb)$, is diagonalisable (because $h^9 = Id$), and has only three eigenvalues: $1$, $e^{2i\pi/3}$, and $e^{-2i\pi/3}$. The dimensions of the corresponding eigenspaces in $H^1(F,\Cb)$ are respectively $8$, $2$ and $2$. In section \ref{example}, we prove the following:
\begin{thm}
The first resonance variety $R_1^1(F)$ of $F$ does not coincide with the tangent cone at $1$, $TC_1V_1^1(F)$, of $V_1^1(F)$. In particular, $F$ is not $1$-formal.
\label{pex}
\end{thm}

\section{Components of the characteristic varieties.}

\subsection{} Let $\Ar = \{H_1, \dots , H_d\}$ be a complex line arrangement in $\Pb^2$. For $i \in \{1,\dots, d\}$, we denote by $f_i$ a linear form such that $H_i = (f_i=0)$. Let $f=f_1 \cdot \ldots \cdot f_d \in \Cb[x,y,z]$ be the homogeneous polynomial associated to $\Ar$. The map $f$ also defines a hypersurface $(f=0)$ of $\Pb^2$. We set $M = \Pb^2 \setminus (f=0)$, and denote by $F$ the Milnor fiber
$$ F = f^{-1}(1) \subset \Cb^3. $$
The natural projection $\pi : \Cb^3\setminus \{0\} \rightarrow \Pb^2$ induces
$$\pi : F \longrightarrow M, $$
which is a $d$-fold covering. The deck transformation group of $F \rightarrow M$ is isomorphic to $\Zb / d\Zb$. We fix a base point $b \in M$, a base point $\tilde{b} \in F$ with $\pi(\tilde{b})=b$, and denote by 
$$ R : \pi_1(M,b) \longrightarrow \pi_1(M,b) / \pi_{\sharp}(\pi_1(F,\tilde{b})) \simeq \Zb / d\Zb $$
the morphism of passage to the quotient.\\
Let us denote by $\gamma_1, \dots, \gamma_d$ elementary loops around, respectively, $H_1,\dots, H_d$. It is known that $\gamma_1, \dots, \gamma_d$ generate $\pi_1(M,b)$. Up to orientation of the $\gamma_i$'s, we have for all $i$, $R(\gamma_i)=1$ (see \cite{CS95}, section $1$).

\subsection{} Let us denote simply, for any complex variety $X$, by $V_1(X)$ the first caracteristic variety $V_1^1(X)$, and similarly $R_1(X)=R_1^1(X)$. According to Arapura \cite{A97}, any positive dimensional irreducible component of $V_1(X)$ is of the form $\rho f^{\ast}H^1(C,\Cb^{\ast})$, where $\rho \in H^1(X,\Cb^{\ast})$ is a torsion character, and $f:X \rightarrow C$ is an admissible map. Moreover, any admissible map $f: X \rightarrow C$ to a curve of negative Euler characteristic (Arapura calls it ``of general type'') gives rise to a component $f^{\ast}H^1(C,\Cb^{\ast})$ of $V_1(X)$.\\
Since $M$ is a line arrangement complement in $\Pb^2$, several kinds of components of $V_1(M)$ can occur, among which:
\begin{itemize}
\item the local components. Fix $p \in \Pb^2$ belonging to $k$ lines of $\Ar$, with $k \ge 3$. Consider the map $M \rightarrow \Pb^1$ which associates to a point $p'$ of $M$ the class of the line passing through $p$ and $p'$. The image of this map is $\Pb^1$ deprived of $k$ points, which has negative Euler characteristic, so, restricted to its image, it is an admissible map of general type, and gives rise to a component of $V_1(M)$.
\item the global components. Recall from \cite{D92}, Chapter 4, that
$$H^1(M,\Cb^{\ast}) \simeq \{(\lambda_1,\dots, \lambda_d)\in (\Cb^{\ast})^d,\, \lambda_1 \cdot \ldots \cdot \lambda_d = 1\}. $$
A global component is a component which contains the origin $1$, and which is included in none of the coordinate hypersurfaces $(\lambda_i=1)$.
\end{itemize}

\subsection{} Let $W$ be a global irreducible component of $V_1(M)$. Then, according to \cite{FY07}, there exists a pencil
$$ \phi : M \longrightarrow S_{\phi} = \Pb^1 \setminus \{a_1, \dots, a_k\} $$
which is an admissible map in the sense of Arapura \cite{A97}, associated to the component $W$. We have $k \ge 3$, so that $\chi(S_{\phi}) < 0$. Actually, $\phi$ extends to a pencil $\hat{\phi}: \Pb^2 \setminus \Br \rightarrow \Pb^1$, where $\Br$ is the base locus of the pencil. There exist exponents $k_1, \dots, k_d$, and indices $m_1<\dots < m_k$ such that
$$ \begin{array}{l}
\hat{\phi}^{-1}(a_1) = (f_1^{k_1}\cdots f_{m_1}^{k_{m_1}} = 0) \\
\hat{\phi}^{-1}(a_2) = (f_{m_1+1}^{k_{m_1+1}} \cdots f_{m_2}^{k_{m_2}} =0) \\
\vdots \\
\hat{\phi}^{-1}(a_k) = (f_{m_{k-1}+1}^{k_{m_{k-1}+1}} \cdots f_d^{k_d} = 0). \end{array} $$

For $i \in \{1,\dots, k\}$, we set $C_i = \phi^{-1}(a_i)$, and $Q_i = \prod_{j=m_{i-1}+1}^{m_i} f_j^{k_j}$. Since $\phi$ is a pencil, the space $<Q_1,\dots, Q_k>$ is of dimension $2$, generated by any two of the polynomials $Q_1$, \dots, $Q_k$. We can suppose that $a_1=[0:1]$ and $a_2 = [1:0]$, so that $\phi$ is the pencil $[Q_1:Q_2]$.\\\\
Since $\phi$ is admissible, its generic fiber is connected. Let $N_0$ be a Zariski-closed subset of $S_{\phi}$ such that for $x \in S_{\phi} \setminus N_0$, $\phi^{-1}(x)$ is connected. Given such a generic $x$, we denote by $F_{\phi}$ a generic fiber of $\phi$.\\

\begin{prop} Suppose that for all $i$, $k_i = 1$, and that there is a point $p$ in the base locus $\Br$ of the pencil $\phi$ such that the multiplicity of $C_j$ in $p$ is $1$ for all $j$.\\
Then there exists a non-constant morphism $\tilde{\phi} : F \rightarrow \tilde{S_{\phi}}$, with connected generic fiber, such that the following diagram is commutative:
$$ \xymatrix{
F \ar^{\tilde{\phi}}[r] \ar_{\pi}[d] & \tilde{S_{\phi}} \ar^{\tilde{\pi}}[d] \\ M \ar^{\phi}[r] & S_{\phi} } $$
Here $\tilde{\pi}$ is induced by the projection $\Cb^2\setminus \{0\} \rightarrow \Pb^1$, and $\tilde{S_{\phi}}$ is a Milnor fiber of a degree $k$ homogeneous polynomial in two variables.
\label{comp}
\end{prop}
\begin{proof} 
We set, for $j \ge 3$, $Q_j = \alpha_j Q_1 + \beta_j Q_2$, and
$$ \tilde{S_{\phi}} = \{ (u,v) \in \Cb^2,\, uv\prod_{j=3}^k (\alpha_ju + \beta_jv) = 1 \}. $$
Then it is straightforward that
$$ \tilde{\phi} : (x,y,z) \longmapsto (Q_1(x,y,z),Q_2(x,y,z)) $$
maps $F$ to $\tilde{S_{\phi}}$, because, since for all $i$, $k_i=1$, $Q_1 \cdots Q_k = f_1 \cdots f_d$. So we have the desired commutative diagram.\\
The image of $\tilde{\phi}$ is a constructible subset of the curve $\tilde{S_{\phi}}$, so equals all $\tilde{S_{\phi}}$ deprived of a finite number of points, whose set of images through $\tilde{\pi}$ we denote by $N_1$. We set $N = N_0 \cup N_1$. It is a finite subset of $S_{\phi}$. We choose a point $y$ in $\tilde{S_{\phi}} \setminus \tilde{\pi}^{-1}(N)$, and show that the fiber $\tilde{\phi}^{-1}(y)$ is connected. We set $F_{\phi} = \phi^{-1}(\tilde{\pi}(y))$, and $\tilde{F_{\phi}} = \pi^{-1}(F_{\phi})$.\\\\
As $F \rightarrow M$ is a $d$-fold covering, $\tilde{F_{\phi}} \rightarrow F_{\phi}$ is a $d$-fold covering. We want to know the number of connected components of $\tilde{F_{\phi}}$, so we study the action of $\pi_1(F_{\phi},p_0)$, with $p_0 \in F_{\phi}$, on $\pi^{-1}(p_0)$. This action is given by the composition of $\pi_1(F_{\phi}) \rightarrow \pi_1(M)$ and $R$.\\
Let $p$ be a point in the base locus of the pencil $\phi$ such that the multiplicity of $C_j$ in $p$ is $1$ for all $j$. Then the closure (in $\Pb^2$) $\overline{F_{\phi}}$ is smooth at $p$. Indeed, $\overline{F_{\phi}}$ is the zero locus of $\alpha Q_1 + \beta Q_2$, for some non-zero $\alpha$ and $\beta$, so it is sufficient that at least one of the two curves $C_1$ or $C_2$ be smooth at $p$, and actually the hypotheses imply both are smooth.\\
Let $\beta_p$ be a small elementary loop in $F_{\phi}$ around $p$. Since $C_1 \cap C_2 = C_1 \cap C_2 \cap \dots \cap C_k$, and with the multiplicity hypotheses at $p$, we can say that $p$ is the intersection of $k$ lines of $\Ar$. Suppose their indices are $i_1,\dots, i_k$. We have
$$ [\beta_p] = [\gamma_{i_1}] + \dots + [\gamma_{i_k}], $$
where $[\gamma]$ denotes the class of the loop $\gamma$ in $H_1(M)$, see \cite{D04}, p.209. So $R(\beta_p) = k$, and $\beta_p$ acts as an element of order $d/k$, which induces $k$ orbits on the fiber $\pi^{-1}(p_0)$. Consequently, the action of $\pi_1(F_{\phi})$ has at most $k$ orbits. This implies $\tilde{F_{\phi}}$ has at most $k$ connected components.\\\\
Let $x \in S_{\phi} \setminus N_1$. The fiber $\tilde{\pi}^{-1}(x)$ is a finite set of exactly $k$ elements, each of which is in the image of $\tilde{\phi}$, so $\tilde{\phi}^{-1}\tilde{\pi}^{-1}(x)$ has at least $k$ connected components.\\\\
Eventually, take $y \in \tilde{S_{\phi}}$ to be a generic point as specified before. Suppose $\tilde{\phi}^{-1}(y)$ (which is not empty) is not connected. Set $x = \tilde{\pi}(y)$. The fiber $\tilde{\phi}^{-1}\tilde{\pi}^{-1}(x)$ has strictly more than $k$ connected components. But
$$ \begin{array}{rcl}
\tilde{\phi}^{-1}\tilde{\pi}^{-1}(x) & = & \pi^{-1}\phi^{-1}(x) \\ & = & \pi^{-1}(F_{\phi}) \\ & = & \tilde{F_{\phi}}, \end{array} $$
which has at most $k$ components. We obtain a contradiction, which means the generic fiber $\tilde{\phi}^{-1}(y)$ is connected. \end{proof}

\subsection{} There is a similar result for local components. Let $W$ be a local component, associated to a point $p$ of multiplicity $k \ge 3$. Let $\phi : M \rightarrow \Pb^1 \setminus \{a_1,\dots, a_k\}$ be the associated admissible map. Let $\tilde{\phi} = \pi \circ \phi$. The diagram is the following:
$$ \xymatrix{
F \ar^{\tilde{\phi}}[rd] \ar_{\pi}[d] &  \\ M \ar^{\phi}[r] & S_{\phi} } $$
\begin{prop} Suppose the arrangement $\Ar$ is not central (not all the lines meet in one point). Then the generic fiber of $\tilde{\phi}: F \rightarrow \Pb^1\setminus \{a_1, \dots, a_k\}$ is connected.
\end{prop}
\begin{proof} Let $i\in \{1, \dots, d\}$ be such that $p \not \in H_i$. Let $a \in \Pb^1\setminus \{a_1, \dots, a_k\}$. Let $p'$ denote the intersection point of $\phi^{-1}(a)$ and $H_i$. Let $\beta$ be a small elementary loop around $p'$ in $\phi^{-1}(a)\setminus \{p'\}$. Then $R(\beta) = R(\gamma_i) = 1$. So this $\beta$ as an element of $\pi_1(\phi^{-1}(a))$ acts with order $d$ on $\tilde{\phi}^{-1}(a)$, hence the action of $\pi_1(\phi^{-1}(a))$ on $\tilde{\phi}^{-1}(a)$ is transitive, and $\tilde{\phi}^{-1}(a)$ is connected. \end{proof}

\noindent{\bf Remark.} Suppose the arrangement is central, and set, for $i \ge 3$, $f_i = \lambda_i f_1 + \mu_i f_2$. Let $a=(\alpha:\beta) \in \Pb^1\setminus \{a_1, \dots, a_k\}$ with $\alpha \neq 0$. Then $\phi$ is the pencil $(x:y:z) \mapsto (f_1(x,y,z):f_2(x,y,z))$, and $\tilde{\phi}^{-1}(a)$ is the set
$$ \{ (x,y,z) \in \Cb^3,\, f_2(x,y,z) = \frac{\beta}{\alpha}f_1(x,y,z) \textrm{ and } f_1f_2\prod_i (\lambda_i f_1 + \mu_if_2)(x,y,z) = 1 \}. $$
In this set, the second condition reduces to 
$$ f_1^d(x,y,z) = \frac{1}{\frac{\beta}{\alpha} \prod_i(\lambda_i + \mu_i \frac{\beta}{\alpha})}. $$
So the fiber $\tilde{\phi}^{-1}(a)$ is a union of $d$ parallel lines in $\Cb^3$, and is not connected.

\section{Example} \label{example}
Consider the line arrangement $\Ar$ of $\Pb^2$, defined in the introduction, given by $(f=0)$, with
$$ f = (x^3-y^3)(x^3-z^3)(y^3-z^3). $$
This section is devoted to the proof of Theorem \ref{pex}, so we suppose from now on $R_1(F) = TC_1V_1(F)$ and derive a contradiction.\\\\
The pencil
$$ \phi : \begin{array}{rcl} M & \longrightarrow & S:=\Pb^1 \setminus \{ [1:0], [0:1], [1:-1] \} \\ (x:y:z) & \longmapsto & (x^3-y^3 : y^3-z^3) \end{array} $$
is an admissible map describing a global irreducible component of the characteristic variety $V_1(M)$. By Proposition \ref{comp}, it comes with a commutative diagram
$$ \xymatrix{ F \ar_{\pi}[d] \ar^{\tilde{\phi}}[r] & \tilde{S} \ar^{\tilde{\pi}}[d] \\ M \ar^{\phi}[r] & S } $$
where $\pi$ and $\tilde{\pi}$ are induced by the natural projections, and
$$ \tilde{S} = \{ (u,v) \in \Cb^2,\, uv(u+v) = 1 \} . $$
Let $E = \tilde{\phi}^{-1} H^1(\tilde{S},\Cb)$. According to Arapura \cite{A97}, since $\chi(\tilde{S}) = -3 < 0$, it is an irreducible component of the tangent cone $TC_1V_1(M)$.
\begin{lemme} There is a mixed Hodge structure on $E$, decomposing as
$$ E = E^{1,1} \oplus E^{1,0} \oplus E^{0,1}, $$
with
$$ dim E^{1,1}=2,\quad dim E^{1,0}=1,\quad dim E^{0,1} = 1. $$
\end{lemme}
\begin{proof} There is a mixed Hodge structure on $H^1(\tilde{S},\Cb)$. Let us refer to \cite{D92}, Theorem C24, to state that the weight filtration for this MHS satisfies $W_0=0$ and $W_2=H^1(\tilde{S},\Cb)$. The Hodge filtration $F$ satisfies $F^0 = H^1(\tilde{S},\Cb)$ and $F^2 =0$. Denote by $H^{1,1}$, $H^{1,0}$, and $H^{0,1}$, respectively, the spaces $F^1 \cap \overline{F}^1$, $W_1 \cap F^1$, and $W_1 \cap \overline{F}^1$. A compactification of $\tilde{S}$ is given by the curve $(uv(u+v)-t^3=0)$ in $\Pb^2$, whose genus is $1$, so point $(iii)$ of Theorem C24 yields $dim W_1 = 2$. Example C26 in \cite{D92} gives
$$ dim H^{1,1} = dim M(uv(u+v))_{1} = 2, $$
where
$$ M(uv(u+v)) = \Cb[u,v] / \left( u^2+2uv, v^2+2uv \right) $$
is the associated graded Milnor algebra.\\
The map $\tilde{\phi}$ is dominant, so induces a monomorphism at the $H^1$ level. Setting for all $i,j$, $E^{i,j} = \tilde{\phi}^{-1} H^{i,j}$, the lemma is proved. \end{proof}
The fiber $F$ is also endowed with a mixed Hodge structure. Since Deligne's construction of a MHS on an algebraic variety is functorial, $\tilde{\phi}$ is a morphism of MHS, and the MHS on $E$ from the lemma is induced by the MHS on $F$. We denote the latter by
$$ H^1(F,\Cb) = H^{1,1}(F) \oplus \underbrace{H^{1,0}(F) \oplus H^{0,1}(F)}_{W_1(F)}. $$
Recall from \cite{DP09}, Theorem 4.1, that
$$ H^{1,1}(F) = H^1(F,\Cb)_1 \simeq H^1(M,\Cb), $$
the middle space denoting the eigenspace associated to the eigenvalue $1$ under the action of monodromy $h^{\ast}$. Also,
$$ W_1(F) = H^1(F,\Cb)_{\neq 1}, $$
with a similar notation. In this example, the only two eigenvalues different from $1$ are $\lambda = e^{2i\pi/3}$ and $\overline{\lambda}$, and $dim W_1(F) = 4$.
\begin{lemme} For all $\alpha \in H^1(F,\Cb)_{\lambda}$ and $\beta \in H^1(F,\Cb)_{\overline{\lambda}}$, $\alpha \cup \beta = 0$. Consequently, $W_1(F) \subset R_1(F)$. \end{lemme}
\begin{proof} Take such $\alpha$ and $\beta$. Then
$$ h^{\ast}(\alpha \cup \beta) = h^{\ast}(\alpha) \cup h^{\ast}(\beta) = \lambda \overline{\lambda} \alpha \cup \beta = \alpha \cup \beta, $$
so $\alpha \cup \beta \in H^2(F,\Cb)_1 \simeq H^2(M,\Cb)$, and this latter space has a pure Hodge structure of weight $4$, whereas $\alpha \cup \beta$ is of weight $2$, so $\alpha \cup \beta = 0$. Now, with the definition of $R_1(F)$, we deduce that $\alpha$ and $\beta$ are in $R_1(F)$, so $H^1(F,\Cb)_{\lambda} \subset R_1(F)$ and $H^1(F,\Cb)_{\overline{\lambda}} \subset R_1(F)$. Let $E'$ be the component of $R_1(F)$ which contains $H^1(F,\Cb)_{\lambda}$. Since we suppose $R_1(F) = TC_1V_1(F)$, $E'$ has a MHS, and then $E' \supset \overline{H^1(F,\Cb)_{\lambda}} = H^1(F,\Cb)_{\overline{\lambda}}$. So $E' \supset W_1(F)$.
\end{proof}
Let $E'$ denote the linear irreducible component of $R_1(F)$ which contains $W_1(F)$. 
Then $E' \cap E \neq \{0\}$, but $E' \cap H^{1,0}(F)$ has dimension $2$ whereas $E \cap H^{1,0}(F) = E^{1,0}$ has dimension $1$, so $E' \neq E$, and this is a contradiction because $E$ is supposed to be a component of $TC_1V_1(F) = R_1(F)$.

\bibliographystyle{amsplain}
\bibliography{biblio.bib}
\end{document}